\documentclass[10pt,reqno]{article}

\usepackage{amssymb}
\usepackage{epsfig}
\usepackage{amsmath}
\usepackage{amsthm}
\usepackage{color}
\definecolor{r}{rgb}{0.9,0.3,0.1}
\definecolor{b}{rgb}{0.1,0.3,0.9}

\newtheorem{theorem}{Theorem}[section]
\newtheorem{lemma}[theorem]{Lemma}

\theoremstyle{remark}
\newtheorem{remark}[theorem]{Remark}

\theoremstyle{definition}
\newtheorem{assumption}[theorem]{Assumption}

\newtheorem{definition}[theorem]{Definition}

\newcommand\cbrk{\text{$]$\kern-.15em$]$}}
\newcommand\opar{\text{\,\raise.2ex\hbox{${\scriptstyle
|}$}\kern-.34em$($}}
\newcommand\cpar{\text{$)$\kern-.34em\raise.2ex\hbox{${\scriptstyle |}$}}\,}

\newcommand{\ga}{\gamma}

\newcommand{\om}{\omega}

\newcommand{\de}{\delta}

\newcommand\bL{\mathbb{L}}
\newcommand\bR{\mathbb{R}}
\newcommand\bH{\mathbb{H}}
\newcommand\bZ{\mathbb{Z}}

\newcommand\cH{\mathcal{H}}

\newcommand\cO{\mathcal{O}}

\newcommand\frH{\mathfrak{H}}

\newcommand{\mysection}[1]{\section{#1}
\setcounter{equation}{0}}

\begin{document}
\setlength{\baselineskip}{16pt}

\title
{A $W^n_2$-Theory of  Elliptic and Parabolic  Partial Differential
Systems in $C^1$ domains}

\author{Kijung Lee\footnote{Department of Mathematics, Ajou University,
Suwon, South Korea 443-749, \,\, kijung@ajou.ac.kr.}
\qquad \hbox{\rm and} \qquad Kyeong-Hun Kim\footnote{Department of
Mathematics, Korea University, 1 Anam-dong, Sungbuk-gu, Seoul, South
Korea 136-701, \,\, kyeonghun@korea.ac.kr. The research of this
author is supported by the Korean Research Foundation Grant funded
by the Korean Government 20090087117}}

\date{}


\maketitle

\begin{abstract}
In this paper second-order  elliptic and parabolic partial
differential systems are considered on $C^1$ domains.  Existence and uniqueness results
are obtained in terms of Sobolev spaces with weights so that we
allow the derivatives of the solutions to blow up near the boundary.
The coefficients of the systems are allowed to substantially
oscillate or blow up near the boundary.

\vspace*{.125in}


\vspace*{.125in}

\end{abstract}



\mysection{Introduction}

In this article we are dealing with the Sobolev space theory of
second-order parabolic and elliptic systems :
\begin{equation}
                       \label{eqn main system}
u^k_t=a^{ij}_{kr}u^r_{x^ix^j}+b^i_{kr}u^r_{x^i}+c_{kr}u^r+f^k,  \quad t>0, x\in \cO
\end{equation}
\begin{equation}
                      \label{elliptic system}
a^{ij}_{kr}u^r_{x^ix^j}+b^{i}_{kr}u^r_{x^i}+c_{kr}u^r+f^k=0, \quad x\in \cO,
\end{equation}
where $\cO$ is a $C^1$ domain in $\bR^d$, $i,j=1,2,\ldots,d$ and
$k,r=1,2,\ldots,d_1$. We used summation notation on repeated indices
$i,j,r$. 

Since the boundary is not supposed to be regular enough, we have to
look for solutions in function spaces with weights, allowing
derivatives of our solutions to blow up near the boundary. In the
framework of H\"older space  such setting leads to investigating
so-called intermediate (or interior) Schauder estimates, which
originated in \cite{DN}.  For results about these estimates the
reader is referred to \cite{DN}, \cite{GH}, \cite{GT} (elliptic
case) and \cite{Fr}, \cite{Li} (parabolic case).

Various Sobolev spaces with weights
 and their
applications to partial differential equations have been
investigated since long ago; we do not want even to  try to collect
all relevant references (some of them can be found in \cite{Ch}).
The reader can find a part of references related to the subject of
this article in   the papers \cite{KK2}, \cite{kr99}
 and \cite{Lo2}, the results of which are extensively
used in what follows.

The main source of our interest in the Sobolev space theory of systems (\ref{eqn main system}) and (\ref{elliptic system})
comes from   \cite{KK2},  \cite{kr99}, \cite{KL2} and \cite{Lo2}, where weighted Sobolev space theory is constructed
for single equations. The goal of this article is to extend   the results for
{\bf{single equations}} in  \cite{KK2}, \cite{kr99}, \cite{KL2} and \cite{Lo2} to the
case of {\bf{the systems}}. We
 prove the uniqueness and existence results   of  systems (\ref{eqn main system}) and (\ref{elliptic system})  in weighted Sobolev spaces under minimal regularity conditions on the coefficients.
 As in the  articles referred above,  our coefficients $a^{ij}_{kr}$
are allowed to substantially oscillate near the boundary, and the
coefficients $b^{i}_{kr},c_{kr}$ are allowed to be unbounded and blow up
near the boundary. For instance, if $d=d_1=1$ and $\cO=(0,\infty)$, then we allow $a:=a^{11}_{11}$
to behave near $x=0$ like $2+\cos|\ln x|^{\alpha}$, where $\alpha\in
(0,1)$ (see Remark \ref{05.18.01}).

However, unlike in those articles, we were able to obtain only
$L_2$-estimates, instead of $L_p$-estimates. This is due to the
difficulty caused by considering systems instead of  single
equations.  For $L_p$-theory, $p>2$, one must overcome tremendous mathematical difficulties rising in the
general settings; one of the main difficulties in the case $p>2$ is
that the arguments we are using in the proofs of  Lemma \ref{a priori 1} and Lemma \ref{a priori 2} below are not
working when $p>2$ since in this case we get extra terms which we
simply can not control.

The organization of the article is as follows. Section \ref{Cauchy}
handles the Cauchy problem. In section \ref{main section} we present
our main results, Theorem \ref{main theorem on domain} and Theorem \ref{theorem elliptic}. In section \ref{section auxiliary}
we develop some auxiliary results.  Theorem \ref{main theorem on domain} and Theorem \ref{theorem elliptic} are proved in section \ref{section 5} and section \ref{section 6}, respectively.

 As usual $\bR^{d}$
stands for the Euclidean space of points $x=(x^{1},...,x^{d})$,
$B_{r}(x)=\{y\in\bR^{d}:|x-y|<r\}$, $B_{r}=B_{r}(0)$,
$\bR^{d}_{+}=\{x\in\bR^{d}:x^{1}>0\}$.
For $i=1,...,d$, multi-indices $\alpha=(\alpha_{1},...,\alpha_{d})$,
$\alpha_{i}\in\{0,1,2,...\}$, and functions $u(x)$ we set
$$
u_{x^{i}}=\frac{\partial u}{\partial x^{i}}=D_{i}u,\quad
D^{\alpha}u=D_{1}^{\alpha_{1}}\cdot...\cdot D^{\alpha_{d}}_{d}u,
\quad|\alpha|=\alpha_{1}+...+\alpha_{d}.
$$
If we write  $c=c(\cdots)$, this means that the constant $c$ depends
only on what are in parenthesis.

\mysection{The system on $\bR^d$}\label{Cauchy}

First we introduce  some solvability results of linear systems
defined on $\bR^d$. These results will be used later for systems
defined on the half space and  bounded $C^1$ domains.

 Let
$C^{\infty}_0=C^{\infty}_0(\mathbb{R}^d;\mathbb{R}^{d_1})$ denote
the set of all $\mathbb{R}^{d_1}$-valued infinitely differentiable
functions with compact support in $\mathbb{R}^d$. By $\mathcal{D}$
we denote the space of $\mathbb{R}^{d}$-valued distributions on
$C^{\infty}_0$; precisely, for $u\in \mathcal{D}$ and $\phi\in
C^{\infty}_0$ we define $(u,\phi)\in \mathbb{R}^{d}$ with components
$(u,\phi)^k=(u^k,\phi^k)$, $k=1,\ldots,d_1$. Each $u^k$ is a usual
$\mathbb{R}$-valued distribution defined on
$C^{\infty}(\mathbb{R}^d;\mathbb{R})$.

We define $L_p=L_p(\mathbb{R}^d;\mathbb{R}^{d_1})$ as the space of
all $\mathbb{R}^{d_1}$-valued functions $u=(u^1,\ldots,u^{d_1})$
satisfying
\[
\|u\|^p_{L_p}:=\sum^{d_1}_{k=1}\|u^k\|^p_{L_p}<\infty.
\]
Let $p \in[2,\infty)$ and $\gamma\in(-\infty,\infty)$. We define the
space of Bessel potential
$H^{\gamma}_p=H^{\gamma}_p(\mathbb{R}^d;\mathbb{R}^{d_1})$ as the
space of all distributions $u$ such that $(1-\Delta)^{n/2}u\in L_p$
where we define each component by
\[
((1-\Delta)^{\gamma/2}u)^k=(1-\Delta)^{\ga/2}u^k
\]
and the norm is given by
\[
\|u\|_{H^{\gamma}_p}:=\|(1-\Delta)^{\ga/2}u\|_{L_p}.
\]
Then, $H^{\gamma}_p$ is a Banach space with the given norm and
$C^{\infty}_0$ is dense in $H^{\gamma}_p$. Note that $H^{\ga}_p$ are
usual Sobolev spaces for $\ga=0,1,2,\ldots$. It is well known that
the first order differentiation operators,
$\partial_i:H^{\gamma}_{p}(\mathbb{R}^d;\mathbb{R})\to
H^{\gamma-1}_p(\mathbb{R}^d;\mathbb{R})$ given by $u\to u_{x^i}$
$(i=1,2,\ldots,d)$, are bounded. On the other hand, for $u\in
H^{\gamma}_{p}(\mathbb{R}^d;\mathbb{R})$, if $\text{supp}\, (u)
\subset (a,b)\times \mathbb{R}^{d-1}$ with $-\infty<a<b<\infty$,
we have
\begin{equation}
                                        \label{eqn 5.1.1}
\|u\|_{H^{\gamma}_{p}(\mathbb{R}^d;\mathbb{R})}\leq
c(d,a,b)\|u_{x^1}\|_{H^{\gamma-1}_{p}(\mathbb{R}^d;\mathbb{R})}
\end{equation}
(see, for instance, Remark 1.13 in \cite{kr99}).

For a fixed time $T$, we define
$$
\bH^{\ga}_p(T):=L_p((0,T],H^{\ga}_p), \quad
\mathbb{L}_p(T):=\bH^0_p(T)
$$
with the norm given by
\[
\|u\|^p_{\mathbb{H}^{\ga}_p(T)}=\int^{T}_0\|u(t)\|^p_{H^{\gamma}_p}dt.
\]
Finally, we set
$U^{\gamma}_p=H^{\gamma-2/p}_p$.

\begin{definition}\label{md}
For a $\mathcal{D}$-valued function $u\in\mathbb{H}^{\ga+2}_p(T)$,
we write $u\in\mathcal{H}^{\ga+2}_p(T)$ if
$u\in\mathbb{H}^{\ga+2}_p(T)$, $u(0,\cdot)\in U^{\gamma+2}_p$
and there exists $f\in \mathbb{H}^{\ga}_p(T)$ such that, for any $\phi\in
C^{\infty}_0$,   the equality
\begin{equation}\label{e}
(u(t,\cdot),\phi)= (u(0,\cdot),\phi)+ \int^t_0( f(s,\cdot),\phi)ds
\end{equation}
holds for all $t\leq T$. In this case, we say that $u_t=f$ \emph{in the sense of distributions}.

The norm in
$\cH^{\gamma+2}_{p}(T)$ is defined by
\[
\|u\|_{\mathcal{H}^{\ga+2}_p(T)}= \|u\|_{\mathbb{H}^{\ga+2}_p(T)}+\|u_t\|_{\mathbb{H}^{\ga}_p(T)}+
\|u(0)\|_{U^{\ga+2}_p}.
\]
\end{definition}

For any $d_1\times d_1$ matrix $C=(c_{kr})$ we let
$$
|C|:=\sqrt{\sum_{k,r}(c_{kr})^2}.
$$
Set $A^{ij}=(a^{ij}_{kr})$.
 Throughout the
article we assume the following.

\begin{assumption}
                     \label{main assumptions}
There exist  constants $\delta,K^j, L>0$ so that

 (i)
\begin{equation}\label{assumption 1}
\de|\xi|^2\le\xi^*_{i} A^{ij}\xi_{j}
\end{equation}
holds for any $t,x$, where $\xi$ is any
(real) $d_1\times d$ matrix, $\xi_i$ is the $i$th column of $\xi$,
and again the summations on $i,j$ are understood.

(ii)
\begin{equation}\label{assumption 2}
\left|A^{1j}\right|\le K^j,  \quad
j=1,2,\ldots,d.
\end{equation}
\end{assumption}

Before we study system (\ref{eqn main system}),  we consider the following system of  equations with constant coefficients:
\begin{equation}
                       \label{eqn system}
u^k_t=a^{ij}_{kr}u^r_{x^ix^j}+f^k,
\quad u^k(0)=u^k_0,
\end{equation}
where $i,j=1,2,\cdots,d$ and $k,r=1,2,\cdots,d_1$;
recall that we are using summation notation on $i,j,r$.

The following $L_2$-theory (even $L_p$-theory) is not new and can be found, for instance, in \cite{Lee}.  However we give a short and  independent  proof for the sake of completeness.

\begin{theorem}
                      \label{thm 1}
Let $a^{ij}_{kr}=a^{ij}_{kr}(t)$, independent of $x$. Then for any $f\in
\bH^{\gamma}_2(T)$ and $u_0\in
U^{\gamma+2}_2$, system (\ref{eqn system}) has a unique solution
$u\in \mathcal{H}^{\gamma+2}_2(T)$, and for this solution
\begin{equation}
                        \label{e 6.5.2}
\|u_{xx}\|_{\bH^{\gamma}_{2}(T)}\leq
c\|f\|_{\bH^{\gamma}_2(T)}+c\|u_0\|_{U^{\gamma+2}_2},
\end{equation}
\begin{equation}
                        \label{e 6.5.3}
\|u\|_{\bH^{\gamma+2}_{2}(T)}\leq
ce^{cT}(\|f\|_{\bH^{\gamma}_2(T)}+\|u_0\|_{U^{\gamma+2}_2}),
\end{equation}
where $c=c(d,d_1,\gamma,\delta,K^j)$.
\end{theorem}

\begin{proof}
By  Theorem 5.1 in \cite{Kr99}, for each $k$, the
equation
$$
u^k_t=\delta \Delta u^k+f^k, \quad u^k(0)=u^k_0
$$
has a solution $u^k\in \mathcal{H}^{\gamma+2}_{2}(T)$. For
$\lambda\in [0,1]$ define
 $A^{ij}_{\lambda}:=(1-\lambda)A^{ij}+\delta_{ij}\lambda\delta I$. Then
$$
|A^{ij}_{\lambda}|\le |A^{ij}|,\quad \delta|\xi|^2\leq
\sum_{i,j}\xi^*A^{ij}_{\lambda}\xi_j
$$
with any $d_1\times d$-matrix $\xi$. Thus having the method of
continuity in mind, we only prove that (\ref{e 6.5.2}) and  (\ref{e
6.5.3}) hold given that a solution $u$ already exists.

 {\bf{Step 1}}.
Assume $\gamma=0$. Applying the chain rule
$d|u^k|^2=2u^kdu^k$ for each $k$,
\begin{equation}
|u^k(t)|^2=|u^k_0|^2+\int^t_0
2u^k(a^{ij}_{kr}u^r_{x^ix^j}+f^k)\,ds,\quad
t>0.\label{square}
\end{equation}

By  integrating with respect to $x$ and using
integrating by parts,
\begin{eqnarray}
&&\int_{\bR^d}|u(t)|^2dx+2\int^t_0\int_{\bR^d}\sum_{i,j}(u_{x^i})^*A^{ij}u_{x^j}dxds\nonumber\\
&=&\int_{\bR^d}|u_0|^2dx
+\int^t_0\int_{\bR^d}2u^*fdxds.
\label{eqn 7.8}
\end{eqnarray}
 Hence, it  follows that
\begin{eqnarray}
&&\int_{\bR^d}|u(t)|^2dx+2\delta\;
\int^t_0\int_{\bR^d}|u_x|^2dxds\nonumber\\
&\leq&\int_{\bR^d}|u_{0}|^2dx + \int^t_0\int_{\bR^d}|u|^2dxds
+\int^t_0\int_{\bR^d}|f|^2dxds.\label{eqn
6.5.5}
\end{eqnarray}
Similarly, for $v=u_{x^n}$ with any $n=1,2,\ldots,d$, we get (see
(\ref{eqn 7.8}))
\begin{eqnarray}
&&\int_{\bR^d}|v(t)|^2dx+2\delta \int^t_0\int_{\bR^d}|v_x|^2dxds\nonumber\\
&\leq&\int_{\bR^d}|(u_0)_{x^n}|^2dx
+\int^t_0\int_{\bR^d}-2v_{x^n}^*f\,dxds\nonumber\\
&\leq& \|u_0\|^2_{U^2_2}+\varepsilon
\|u_{xx}\|^2_{\bL_2(t)}+c\|f\|^2_{\bL_2(t)}.\label{v_x}
\end{eqnarray}
Choosing small $\varepsilon$ and considering all $n$, we have
(\ref{e 6.5.2}). Now, (\ref{v_x}), (\ref{eqn 6.5.5}) and Gronwall's
inequality easily lead to (\ref{e 6.5.3}).

{\bf{Step 2}}. Let $\gamma\neq 0$. The results of this case easily
 follow from the fact that
$(1-\Delta)^{\mu/2}:H^{\gamma}_p\to H^{\gamma-\mu}_p$ is an isometry
for any $\gamma,\mu\in \bR$ when $p\in (1,\infty)$; indeed, $u\in
\cH^{\gamma+2}_2(T)$ is a solution of (\ref{eqn system}) if and only
if $v:=(1-\Delta)^{\gamma/2}u\in \cH^2_2(T)$ is a solution of
(\ref{eqn system}) with $(1-\Delta)^{\gamma/2}f, (1-\Delta)^{\gamma/2}u_0$ in place of $f,
u_0$ respectively. Moreover, for instance,
$$
\|u\|_{\bH^{\gamma+2}_{2}(T)}=\|v\|_{\bH^{2}_{2}(T)}\leq
ce^{cT}\left(\|(1-\Delta)^{\gamma/2}f\|_{\bL_2(T)}
+\|(1-\Delta)^{\gamma/2}u_0\|_{U^2_2}\right)
$$
$$
=ce^{cT}\left(\|f\|_{\bH^{\gamma}_2(T)}+\|u_0\|_{U^{\gamma+2}_2}\right).
$$
The theorem is proved.

\end{proof}

Theorem \ref{thm 1} is extended to the systems with variable
coefficients in the followings.

Fix $\mu>0$. For  $\gamma \in \bR$ define
 $|\gamma|_+=|\gamma|$ if  $|\gamma|=0,1,2,\cdots$ and
 $|\gamma|_+=|\gamma|+\mu$ otherwise. Also define
 $$
 B^{|\gamma|_+}=\begin{cases} B(\bR) &: \quad  \gamma=0\\
 C^{|\gamma|-1,1}(\bR) &: \quad |\gamma|=1,2,... \\
C^{|\gamma|+\mu}(\bR) &: \quad \text{otherwise},
 \end{cases}
$$
where $B$ is the space of bounded functions,  and $C^{|\gamma|-1,1}$
and $C^{|\gamma|+\mu}$ are  usual H\"older spaces.

Consider the system with variable coefficients:
\begin{equation}
                       \label{eqn system2}
u^k_t=a^{ij}_{kr}u^r_{x^ix^j}+b^i_{kr}u^r_{x^i}+c_{kr}u^r+f^k, \quad
u^k(0)=u^k_0.
\end{equation}

\begin{theorem}
                      \label{thm 2}
Assume that the coefficients $a^{ij}_{kr}$ are
uniformly continuous in $x$, that is, for any $\varepsilon>0$ there
exists $\delta=\delta(\varepsilon)>0$ so that for any  $t>0$,
$i,j,k,r$,
$$
|a^{ij}_{kr}(t,x)-a^{ij}_{kr}(t,y)|<\varepsilon,
\quad \text{if}\quad |x-y|<\delta.
$$
Also, assume for any $t>0$, $i,j,k,r$,
$$
|a^{ij}_{kr}(t,\cdot)|_{|\gamma|_+}+|b^i_{kr}(\om,t,\cdot)|_{|\gamma|_+}
+|c_{kr}(\om,t,\cdot)|_{|\gamma|_+}<L.
$$
Then for any $f\in \bH^{\gamma}_2(T)$ and $u_0\in U^{\gamma+2}_2$, system
(\ref{eqn system2}) has a unique solution $u\in
\mathcal{H}^{\gamma+2}_2(T)$, and for this solution we have
$$
\|u\|_{\bH^{\gamma+2}_{2}(T)}\leq
ce^{cT}(\|f\|_{\bH^{\gamma}_2(T)}+\|u_0\|_{U^{\gamma+2}_2}),
$$
where $c=c(d,d_1,\gamma,\delta,K^i,L)$.
\end{theorem}

\begin{proof}
This is an easy extension of Theorem \ref{thm 1} and can be proved
by repeating  the proof of Theorem 5.1 in \cite{Kr99}, where the
theorem is proved when   $d_1=1$. We leave the details to the reader.
\end{proof}

\mysection{The system on  $\mathcal{O} \subset \bR^d$}\label{main
section}

\begin{assumption}
                                         \label{assumption domain}

The domain $\cO$  is of class $C^{1}_{u}$. In other words, for any
$x_0 \in \partial \cO$, there exist constants $r_0, K_0\in(0,\infty)
$ and  a one-to-one continuously differentiable mapping $\Psi$ of
 $B_{r_0}(x_0)$ onto a domain $J\subset\bR^d$ such that

(i) $J_+:=\Psi(B_{r_0}(x_0) \cap \cO) \subset \bR^d_+$ and
$\Psi(x_0)=0$;

(ii)  $\Psi(B_{r_0}(x_0) \cap \partial \cO)= J \cap \{y\in
\bR^d:y^1=0 \}$;

(iii) $\|\Psi\|_{C^{1}(B_{r_0}(x_0))}  \leq K_0 $ and
$|\Psi^{-1}(y_1)-\Psi^{-1}(y_2)| \leq K_0 |y_1 -y_2|$ for any $y_i
\in J$;

(iv)   $\Psi_{x}$ is uniformly continuous in for $B_{r_{0}}(x_{0})$.
\end{assumption}

To proceed further we introduce some well known results from
\cite{GH} and \cite{KK2} (also, see \cite{La} for details).

\begin{lemma}
                                           \label{lemma 10.3.1}
Let the domain $\cO$ be of class $C^{1}_{u}$. Then

(i) there is a bounded real-valued function $\psi$ defined in
$\bar{\cO} $  such that the functions $\psi(x)$ and
$\rho(x):=\text{dist}(x,\partial \cO)$ are comparable  in the part
of
  a neighborhood of $\partial \cO$ lying in $\cO$. In other words, if $\rho(x)$ is
sufficiently small, say $\rho(x)\leq 1$, then $N^{-1}\rho(x) \leq
\psi(x) \leq N\rho(x)$ with some constant
 $N$ independent of $x$,

 (ii) for any  multi-index $\alpha$,
\begin{equation}
                                                             \label{03.04.01}
\sup_{\cO} \psi ^{|\alpha|}(x)|D^{\alpha}\psi_{x}(x)| <\infty.
\end{equation}

\end{lemma}

To describe the assumptions of $f$  we use the Banach spaces
introduced in \cite{KK2} and \cite{Lo2}.
 Let $\zeta\in C^{\infty}_{0}(\bR_{+})$
be a   function satisfying
\begin{equation}
                                       \label{eqn 5.6.5}
\sum_{n=-\infty}^{\infty}\zeta(e^{n+x})>c>0, \quad \forall x\in \bR,
\end{equation}
where $c$ is a constant.  Note that  any nonnegative function
$\zeta$, $\zeta>0$ on $[1,e]$, satisfies (\ref{eqn 5.6.5}). For $x\in \cO$ and $n\in\bZ=\{0,\pm1,...\}$
define
$$
\zeta_{n}(x)=\zeta(e^{n}\psi(x)).
$$
Then  we have $\sum_{n}\zeta_{n}\geq c$ in $\cO$ and
\begin{equation*}
\zeta_n \in C^{\infty}_0(\cO), \quad |D^m \zeta_n(x)|\leq
N(m)e^{mn}.
\end{equation*}
For $\theta,\gamma \in \bR$, let $H^{\gamma}_{p,\theta}(\cO)$ be the
set of all distributions $u=(u^1,u^2,\cdots u^{d_1})$  on $\cO$ such
that
\begin{equation}
                                                 \label{10.10.03}
\|u\|_{H^{\gamma}_{p,\theta}(\cO)}^{p}:= \sum_{n\in\bZ} e^{n\theta}
\|\zeta_{-n}(e^{n} \cdot)u(e^{n} \cdot)\|^p_{H^{\gamma}_p} < \infty.
\end{equation}

It is known (see, for instance, \cite{Lo2}) that up to equivalent
norms the space $H^{\gamma}_{p,\theta}(\cO)$ is independent of the
choice of $\zeta$ and $\psi$. Moreover if $\gamma=n$ is a
non-negative integer then
\begin{equation}
                              \label{eqn 02.09.1}
\|u\|^p_{H^{\gamma}_{p,\theta}(\cO)} \sim \sum_{k=0}^n
\sum_{|\alpha|=k}\int_{\cO} |\psi^kD^{\alpha}u(x)|^p
\psi^{\theta-d}(x) \,dx.
\end{equation}

Denote $\rho(x,y)=\rho(x)\wedge \rho(y)$ and
$\psi(x,y)=\psi(x)\wedge \psi(y)$. For
  $n \in\bZ$, $\mu \in(0,1]$
 and $k=0,1,2,...$, define
$$
|u|_{C}=\sup_{\cO}|u(x)|, \quad [u]_{C^{\mu}}=\sup_{x\neq
y}\frac{|u(x)-u(y)|}{|x-y|^{\mu}}.
$$
\begin{equation}
                           \label{eqn 5.6.2}
[u]^{(n)}_{k}=[u]^{(n)}_{k,\cO} =\sup_{\substack{x\in \cO\\
|\beta|=k}}\psi^{k+n}(x)|D^{\beta}u(x)|,
\end{equation}
\begin{equation}
                              \label{eqn 5.6.3}
[u]^{(n)}_{k+\mu}=[u]^{(n)}_{k+\mu,\cO} =\sup_{\substack{x,y\in \cO
\\ |\beta|=k}}
\psi^{k+\mu+n}(x,y)\frac{|D^{\beta}u(x)-D^{\beta}u(y)|}
{|x-y|^{\mu}},
\end{equation}
$$
|u|^{(n)}_{k}=|u|^{(n)}_{k,\cO}=\sum_{j=0}^{k}[u]^{(n)}_{j,\cO},
\quad |u|^{(n)}_{k+\mu}=
 |u|^{(n)}_{k+\mu,\cO}=|u|^{(n)}_{k, \cO}+
[u]^{(n)}_{k+\mu,\cO}.
$$
In case $\cO=\bR_+$, we also define the norm
$|u|^{(n)*}_{k}=|u|^{(n)*}_{k,\bR_+}$ by using
$\rho(x)(=x^1)$ and $\rho(x)\wedge \rho(y)$ in place of $\psi(x)$
and $\psi(x,y)$ respectively in (\ref{eqn 5.6.2}) and (\ref{eqn
5.6.3}).

Below we collect some other properties of spaces
$H^{\gamma}_{p,\theta}(\cO)$.

\begin{lemma} $(\cite{kr99})$ Let $d-1<\theta<d-1+p$.

                  \label{lemma 1}
(i) Assume that $\gamma-d/p=m+\nu$ for some $m=0,1,\cdots$ and
$\nu\in (0,1]$.  Then for any $u\in H^{\gamma}_{p,\theta}(\cO)$ and $i\in
\{0,1,\cdots,m\}$, we have
$$
|\psi^{i+\theta/p}D^iu|_{C}+[\psi^{m+\nu+\theta/p}D^m
u]_{C^{\nu}}\leq c \|u\|_{ H^{\gamma}_{p,\theta}(\cO)}.
$$

(ii) Let $\alpha\in \bR$, then
$\psi^{\alpha}H^{\gamma}_{p,\theta+\alpha
p}(\cO)=H^{\gamma}_{p,\theta}(\cO)$,
$$
\|u\|_{H^{\gamma}_{p,\theta}(\cO)}\leq c
\|\psi^{-\alpha}u\|_{H^{\gamma}_{p,\theta+\alpha p}(\cO)}\leq
c\|u\|_{H^{\gamma}_{p,\theta}(\cO)}.
$$

(iii) There is a constant $c=c(d,p,\gamma,\theta)$ so that
$$
\|a f\|_{H^{\gamma}_{p,\theta}(\cO)}\leq
c|a|^{(0)}_{|\gamma|_+}|f|_{H^{\gamma}_{p,\theta}(\cO)}.
$$

(iv) $\psi D, D\psi: H^{\gamma}_{p,\theta}(\cO)\to
H^{\gamma-1}_{p,\theta}(\cO)$ are bounded linear operators, and
$$
\|u\|_{H^{\gamma}_{p,\theta}(\cO)}\leq
c\|u\|_{H^{\gamma-1}_{p,\theta}(\cO)}+c \|\psi
Du\|_{H^{\gamma-1}_{p,\theta}(\cO)}\leq c
\|u\|_{H^{\gamma}_{p,\theta}(\cO)},
$$
$$
\|u\|_{H^{\gamma}_{p,\theta}(\cO)}\leq
c\|u\|_{H^{\gamma-1}_{p,\theta}(\cO)}+c \|D\psi
u\|_{H^{\gamma-1}_{p,\theta}(\cO)}\leq c
\|u\|_{H^{\gamma}_{p,\theta}(\cO)}.
$$

\end{lemma}

Denote
$$
\bH^{\gamma}_{p,\theta}(\cO,T)=L_p(
[0,T],H^{\gamma}_{p,\theta}(\cO)), \quad \bL_{p,\theta}(\cO,T)=\bH^{0}_{p,\theta}(\cO,T)
$$
$$
U^{\gamma}_{p,\theta}(\cO)=
\psi^{1-2/p}H^{\gamma-2/p}_{p,\theta}(\cO)).
$$

\begin{definition}
We write   $u\in \frH^{\gamma+2}_{p,\theta}(\cO,T)$  if
 $u=(u^1,\cdots, u^{d_1})\in \psi\bH^{\gamma+2}_{p,\theta}(\cO,T)$,
$u(0,\cdot) \in U^{\gamma+2}_{p,\theta}(\cO)$ and  for some $f \in
\psi^{-1}\bH^{\gamma}_{p,\theta}(\cO,T)$,  it holds that  $u_t=f$
in the sense of distributions. The norm in  $
\frH^{\gamma+2}_{p,\theta}(\cO,T)$ is introduced by
$$
\|u\|_{\frH^{\gamma+2}_{p,\theta}(\cO,T)}=
\|\psi^{-1}u\|_{\bH^{\gamma+2}_{p,\theta}(\cO,T)} + \|\psi
u_t\|_{\bH^{\gamma}_{p,\theta}(\cO,T)}  +
\|u(0,\cdot)\|_{U^{\gamma+2}_{p,\theta}(\cO)}.
$$

\end{definition}

The following result  is due to N.V.Krylov (see \cite{Kr01} and
\cite{Kim04-1}).

\begin{lemma}
                             \label{lemma 15.05}
 Let $p\geq 2$. Then there exists a constant
$c=c(d,p,\theta,\gamma,T)$ such that
$$
 \sup_{t\leq T}\|u(t)\|_{H^{\gamma+1}_{p,\theta}(\cO)}\leq c
\|u\|_{\frH^{\gamma+2}_{p,\theta}(\cO,T)}.
$$
In particular, for any $t\leq T$,
$$
\|u\|^p_{\bH^{\gamma+1}_{p,\theta}(\cO,t)}\leq c \int^t_0
\|u\|^p_{\frH^{\gamma+2}_{p,\theta}(\cO,s)}ds.
$$
\end{lemma}

\begin{assumption}
             \label{assumption regularity}
(i) The functions  $a^{ij}_{kr}(t,\cdot)$
 are  point-wise  continuous in $\cO$, that is, for any $\varepsilon >0, x\in \cO$ there exists $\delta=\delta(\varepsilon,x)$ so that
$$
|a^{ij}_{kr}(t,x)-a^{ij}_{kr}(t,y)|<\varepsilon
$$
whenever $y\in \cO$ and $|x-y|<\delta$.

(ii) There is   control on the behavior of $a^{ij}_{kr}$,
$b^i_{kr}$ and  $c_{kr}$ near
$\partial \cO$, namely,
\begin{equation}
                                                \label{12.10.1}
\lim_{\substack{\rho(x)\to 0\\
x\in \cO}}\sup_{\substack{y\in \cO \\|x-y|\leq\rho(x,y)}} \sup_{t} |a^{ij}_{kr}(t,x)-a^{ij}_{kr}(t,y) | =0.
\end{equation}
\begin{equation}
                                                       \label{05.04.01}
\lim_{\substack{\rho(x)\to0\\
x\in \cO}} \sup_{t}[\rho(x)|b^i_{kr}(t,x)|+\rho^{2}(x)|c_{kr}(t,x)|]=0.
\end{equation}

(iii) For any $t>0$,
$$
|a^{ij}_{kr}(t,\cdot)|^{(0)}_{|\gamma|_+}
+|b^{i}_{kr}(t,\cdot)|^{(1)}_{|\gamma|_+}
+|c_{kr}(t,\cdot)|^{(2)}_{|\gamma|_+}  \leq L.
$$
\end{assumption}

\begin{remark}
          \label{05.18.01}
            It is easy to see
 that (\ref{12.10.1}) is much weaker than uniform continuity condition.
   For instance,
if $\delta\in(0,1)$, $d=d_1=1$, and $\cO=\bR_{+}$, then the function
$a(x)$ equal to $2+\sin (|\ln x|^{\delta})$ for $0<x\leq1/2$
satisfies (\ref{12.10.1}). Indeed, if $x,y>0$ and $ |x-y|\leq
x\wedge y $, then
 $$
|a(x)-a(y)|=|x-y||a'(\xi)|,
$$
where $\xi$ lies between $x$ and $y$. In addition, $|x-y|\leq
x\wedge y \leq\xi\leq2(x\wedge y)$,  and $ \xi|a'(\xi)|\leq
|\ln[2(x\wedge y)]|^{\delta-1}\to0$ as $x\wedge y\to0$.

 Also observe that (\ref{05.04.01}) allows the coefficients
$b^i_{kr}$ and  $c_{kr}$ to blow up near the boundary at a
certain rate.
\end{remark}

\vspace{3mm}
Now, for each $i,j$,
 we define the symmetric part ($S^{ij}$) and the diagonal part ($S^{ij}_d$) of $A^{ij}$ as follows:
$$
S^{ij}=(s^{ij}_{kr}):=(A^{ij}+(A^{ij})^*)/2, \quad \quad S^{ij}_d=(s^{ij}_{d,kr}):=(\delta_{kr}a^{ij}_{kr})=(\delta_{kr}s^{ij}_{kr}).
$$
Also define
$$
H^{ij}:=A^{ij}-(A^{ij})^*, \quad
 S^{ij}_o=S^{ij}-S^{ij}_d.
 $$

 Assume there exist constants $\alpha,\bar{\alpha},\beta_1,\cdots,\beta_d\in [0,\infty)$ so that
\begin{equation}
                         \label{eqn 01.26.1}
  |H^{1j}|\leq \beta^j \quad \forall j=1,2,\ldots,d, \quad\quad |S^{11}_o|\leq \alpha,
\end{equation}
\begin{equation}
                         \label{eqn 01.26.2}
  \xi^*_iS^{ij}_o\xi_j \leq \bar{\alpha}|\xi|^2,
\end{equation}
for any   (real) $d_1\times d$ matrix $\xi$. Here $\xi_i$ is the $i$th column of $\xi$,
and again the summations on $i,j$ are understood. Denote
 $$
 K:=\sqrt{\sum_j (K^j)^2}, \quad \beta=\sqrt{\sum_j (\beta^j)^2}.
 $$

 \begin{assumption}
                             \label{assumption theta}
One of the following four conditions is satisfied:
\begin{equation}
                       \label{theta 11}
 \theta\in \left(d-\frac{\delta}{2K-\delta},\,\,
d+\frac{\delta}{2K+\delta}\right);
\end{equation}
\begin{equation}
             \label{con 11}
\theta\in (d-1,d], \quad
2\delta(d+1-\theta)^2-2(d+1-\theta)(d-\theta)\beta-4(d-\theta)(d+1-\theta)K^1>0;
\end{equation}
\begin{equation}
                   \label{con 44}
\theta\in (d-1,d], \quad (\delta-\bar{\alpha})-\frac{(d-\theta)}{(d+1-\theta)}(2\delta-\beta-2\alpha)>0;
\end{equation}
\begin{equation}
                   \label{con 22}
\theta\in [d,d+1), \quad 8(d+1-\theta)\delta^2-(\theta-d)\beta^2>0.
\end{equation}
 \end{assumption}

\begin{remark}
 (i) If $A^{1j}$ are symmetric, i.e., $\beta=0$, then
(\ref{con 11}) combined with (\ref{con 22}) is $\theta\in
(d-\frac{\delta}{2K^1-\delta},d+1)$ which is weaker than (\ref{theta
11}).

(ii) If $A^{ij}$ are diagonal matrices, that is if $\alpha=\beta^i=0$, then (\ref{con 11}) combined with (\ref{con 22})
is $\theta\in  (d-1,d+1)$. This is the case when the system is not correlated.

(iii)
We also mention that if $\theta\not\in (d-1,d+1)$ then Theorem
\ref{main theorem on domain} is false even for the heat equation
$u_t=\Delta u+f$ (see \cite{kr99}).
\end{remark}

Here are the main results of this article. The proofs of the theorems will be given in section \ref{section 5}
  and section \ref{section 6} after we develop some auxiliary results on $\bR^d_+$ in section \ref{section auxiliary}.

\begin{theorem}
                    \label{main theorem on domain}
Let $\gamma \geq 0$ and  $\cO$ be bounded. Also let Assumptions
\ref{main assumptions}, \ref{assumption domain}, \ref{assumption
regularity} and \ref{assumption theta} hold. Then for any $f\in
\psi^{-1}\bH^{\gamma}_{2,\theta}(\cO,T),\, u_0\in
U^{\gamma+2}_{2,\theta}(\cO)$, system (\ref{eqn system2}) admits a
unique solution $u\in \frH^{\gamma+2}_{2,\theta}(\cO,T)$, and for
this solution
\begin{equation}
                        \label{a priori domain}
\|\psi^{-1}u\|_{\bH^{\gamma+2}_{2,\theta}(\cO,T)}\leq ce^{cT}\left(\|\psi
f\|_{\bH^{\gamma}_{2,\theta}(\cO,T)}
+\|u_0\|_{U^{\gamma+2}_{2,\theta}(\cO)}\right),
\end{equation}
where $c=c(d,\delta,\theta,K,L)$.
\end{theorem}

\begin{theorem}
                           \label{theorem elliptic}
Let $\gamma\geq 0$ and  $\cO$ be bounded. Assume $a^{ij}_{kr},b^{i}_{kr},c_{kr}$ are independent of $t$ and $\lambda^k$ are
 sufficiently large constants (actually, any constants bigger than
$c$ from (\ref{a priori domain})). Under the assumptions of Theorem
\ref{main theorem on domain}, for any $f\in
\psi^{-1}H^{\gamma}_{2,\theta}(\cO)$ there is a unique $u\in \psi
H^{\gamma+2}_{2,\theta}(\cO)$ such that in $\cO$,
\begin{equation}
                     \label{eqn elliptic}
a^{ij}_{kr}u^r_{x^ix^j}+b^{i}_{kr}u^r_{x^i}+c_{kr}u^r-\lambda^ku^k+f^k=0.
\end{equation}
Furthermore,
\begin{equation}
                         \label{a priori elliptic}
\|\psi^{-1}u\|_{H^{\gamma+2}_{2,\theta}(\cO)}\leq N\|\psi
f\|_{H^{\gamma}_{2,\theta}(\cO)},
\end{equation}
where the constant $N$ is independent of $f$.
\end{theorem}

\begin{remark}
Actually  Theorem \ref{main theorem on domain} and Theorem \ref{theorem elliptic} hold even for
$\gamma<0$.  Using  results for the case $\gamma\geq 0$, repeat the arguments in the proof of  Theorem 2.10 in \cite{KK2},
where the theorems are proved when $d_1=1$. We leave the details to the reader. Also
 by inspecting the proofs carefully one can check that  the above two theorems hold true even if $\cO$ is not bounded.

\end{remark}

\mysection{Auxiliary results: some results on $\bR^d_+$}
\label{section auxiliary}

In this section we develop some results for the systems defined on
$\bR^d_+$. Here we use   the Banach spaces $H^{\gamma}_{p,\theta}$,
$\bH^{\gamma}_{p,\theta}(T)$ and $\frH^{\gamma}_{p,\theta}(T)$
 defined on $\bR^d_+$.
  They are defined on the basis of (\ref{10.10.03})
by  formally taking $\psi(x)=x^{1}$, so that $\zeta_{-n}(e^{n}
x)=\zeta(x)$ and
\begin{equation*}
\|u\|_{H^{\gamma}_{p,\theta} }^{p}:= \sum_{n\in\bZ} e^{n\theta}
\|u(e^{n} \cdot)\zeta(\cdot) \|^p_{H^{\gamma}_p} < \infty.
\end{equation*}
Observe that  the spaces $H^{\gamma}_{p,\theta}(\bR^d_+)$ and
$H^{\gamma}_{p,\theta}$ are different since $\psi$ is bounded.
Actually for any nonnegative function $\xi=\xi(x^1)\in
C^{\infty}_0(\bR^1)$ so that  $\xi=1$ near $x^1=0$ we have
\begin{equation}
                             \label{eqn 9.9}
\|u\|_{H^{\gamma}_{p,\theta}(\bR^d_+)}\sim \left(\|\xi
u\|_{H^{\gamma}_{p,\theta}}+\|(1-\xi)u\|_{H^{\gamma}_p}\right).
\end{equation}
Also, it is known (see \cite{kr99}) that  for any $\eta\in
C^{\infty}_0(\bR^d_+)$,
\begin{equation}
                            \label{eqn 5.6.1}
\sum_{n=-\infty}^{\infty}
e^{n\theta}\|u(e^n\cdot)\eta\|^p_{H^{\ga}_p} \leq c
\sum_{n=-\infty}^{\infty}
e^{n\theta}\|u(e^n\cdot)\zeta\|^p_{H^{\ga}_p},
\end{equation}
where $c$ depends only on $d,d_1,\gamma,\theta,p,\eta,\zeta$.
Furthermore,
 if $\gamma=n$ is a nonnegative integer then (see (\ref{eqn 02.09.1}))
\begin{equation}
                              \label{eqn 02.09.2}
\|u\|^p_{H^{\gamma}_{p,\theta}} \sim \sum_{k=0}^n
\sum_{|\alpha|=k}\int_{\bR^d_+} |\psi^kD^{\alpha}u(x)|^p
(x^1)^{\theta-d}(x) \,dx.
\end{equation}

Let $M^{\alpha}$ be the operator of multiplying $(x^1)^{\alpha}$ and
$M=M^1$.

\begin{lemma}
                          \label{collection half}
 The assertions (i)-(iv) in Lemma \ref{lemma 1} hold true if one
formally replaces  $H^{\gamma}_{p,\theta}(\cO)$ and
 $\psi$ by $H^{\gamma}_{p,\theta}$ and $M$, respectively.

\end{lemma}

We need the following three lemmas to prove the main result of this section.

\begin{lemma}
                    \label{lemma 2.05.10}
Let $a^{ij}_{kr}=a^{ij}_{kr}(t)$, independent of $x$.  Assume that  $f\in M^{-1}\bH^{\gamma}_{2,\theta}(T), u(0)\in
U^{\gamma+2}_{2,\theta}$ and $u\in M\bH^{\gamma+1}_{2,\theta}(T)$ is
a solution of  system (\ref{eqn system}) on $[0,T]\times
\mathbb{R}_+$, then $u\in M\bH^{\gamma+2}_{2,\theta}(T)$ and
\begin{equation}
                              \label{eqn 5.1}
\|M^{-1}u\|_{\bH^{\gamma+2}_{2,\theta}(T)}\leq
c\|M^{-1}u\|_{\bH^{\gamma+1}_{2,\theta}(T)}+c\|Mf\|_{\bH^{\gamma}_{2,\theta}(T)}
+c\|u(0)\|_{U^{\gamma+2}_{2,\theta}},
\end{equation}
where $c=c(d,d_1,\gamma,\theta,\delta,K,L)$.
\end{lemma}
\begin{proof}
By Lemma \ref{collection half}  and (\ref{eqn 5.1.1}),
\begin{eqnarray*}
 \|M^{-1}u\|^2_{\bH^{\gamma+2}_{2,\theta}(T)}
&\leq&
c\sum_{n}e^{n(\theta-2)}\|u(t,e^nx)\zeta(x)\|^2_{\bH^{\gamma+2}_2(T)}\\
&=&
c\sum_{n}e^{n\theta}\|u(e^{2n}t,e^nx)\zeta(x)\|^2_{\bH^{\gamma+2}_2(e^{-2n}T)}\\
&\leq&
c\sum_{n}e^{n\theta}\|(u(e^{2n}t,e^nx)\zeta(x))_{x^1x^1}\|^2_{\bH^{\gamma}_2(e^{-2n}T)}.
\end{eqnarray*}
Denote
$$
v_n(t,x)=u(e^{2n}t,e^nx)\zeta(x),\quad
a^{ij}_{n,kr}(t)=a^{ij}_{kr}(e^{2n}t).
$$
Then since $v_n$ has compact support in $\bR^d_+$, $v_n$ is in
$\bH^{\gamma+1}_2(e^{-2n}T)$ and satisfies
$$
(v^k_n)_t= a^{ij}_{n,kr}(v^r_n)_{x^ix^j}+f^k_n,
\quad v^k_n(0,x)=\zeta(x)u^k_0(e^nx),
$$
where
$$
f^k_n=-2e^n a^{ij}_{n,kr} u^r_{x^i}(e^{2n}t,e^nx)\zeta_{x^j}(x) -a^{ij}_{n,kr}
u^r(e^{2n}t,e^nx)\zeta_{x^ix^j}(x)+e^{2n}f^k(e^{2n}t,e^nx)\zeta(x).
$$

By Theorem \ref{thm 1}, $v_n$ is in $\bH^{\gamma+2}_2(e^{-2n}T)$ and
$$
\|(v_{n})_{xx}\|^2_{\bH^{\gamma}_2(e^{-2n}T)}\leq
c(d,d_1,\gamma,\delta,K,L)(\|f_n\|^2_{\bH^{\gamma}_2(e^{-2n}T)}
+\|\zeta(x)u_0(e^nx)\|^2_{U^{\gamma+2}_2}).
$$
 Thus by (\ref{eqn 5.6.1})
and  Lemma \ref{collection half},
\begin{eqnarray*}
&&\sum_{n}e^{n\theta}\|(u(e^{2n}t,e^nx)\zeta(x))_{xx}\|^2_{\bH^{\gamma}_2(e^{-2n}T)}\\
&\leq&
c\sum_{n}e^{n\theta}\|u_{x}(t,e^n\cdot)\zeta_{x}\|^2_{\bH^{\gamma}_2(T)}
+c\sum_{n}e^{n(\theta-2)}\|u(t,e^n\cdot)\zeta_{xx}\|^2_{\bH^{\gamma}_2(T)}\\&&+
c\sum_{n}e^{n(\theta+2)}\|f(t,e^n\cdot)\zeta\|^2_{\bH^{\gamma}_2(T)}
+c\sum_{n}e^{n\theta}\|u_0(t,e^nx)\zeta\|^2_{U^{\gamma+2}_2}\\
&\leq & c \|u_x\|^2_{\bH^{\gamma}_{2,\theta}(T)}+c\|M^{-1}u\|^2_{\bH^{\gamma}_{2,\theta}(T)}
+c\|Mf\|^2_{\bH^{\gamma}_{2,\theta}(T)}+c\|u_0\|^2_{U^{\gamma+2}_{2,\theta}}\\
&\leq &c\|M^{-1}u\|^2_{\bH^{\gamma+1}_{2,\theta}(T)}+c\|M
f\|^2_{\bH^{\gamma}_{2,\theta}(T)}+c\|u_0\|^2_{U^{\gamma+2}_{2,\theta}}.
\end{eqnarray*}
The lemma is proved.
\end{proof}

It follows from the above lemma that if $\gamma \geq 0$, then
$$
\|M^{-1}u\|_{\bH^{\gamma+2}_{2,\theta}(T)}\leq c\|M^{-1}u\|_{\bL_{2,\theta}(T)}
+c\|Mf\|_{\bH^{\gamma}_{2,\theta}(T)}+c\|u_0\|_{U^{\gamma+2}_{2,\theta}}.
$$
Thus to get a priori estimate, we only need to estimate $\|M^{-1}u\|_{\bL_{2,\theta}(T)}$ in terms of $f$ and $u_0$.

\begin{lemma}
                \label{a priori 1}
Let $a^{ij}_{kr}=a^{ij}_{kr}(t)$, independent of $x$. Assume
\begin{equation}
                     \label{theta 1}
\theta\in \left(d-\frac{\delta}{2K-\delta},\,\,
d+\frac{\delta}{2K+\delta}\right)
\end{equation}
and  $u\in M\bH^{1}_{2,\theta}(T)$ is a solution of $($\ref{eqn
system}$)$ so that $u\in C([0,T],C^2_0((1/N,N)\times
\{x':|x'|<N\}))$ for some $N>0$. Then we have
\begin{equation}
               \label{eqn main}
\|M^{-1}u\|^2_{\bL_{2,\theta}(T)}\leq
c_0(\|Mf\|^2_{\bL_{2,\theta}(T)}+\|u_0\|^2_{U^1_{2,\theta}}),
\end{equation}
where $c_0=c_0(d,\delta,\theta,K,L)$.
\end{lemma}

\begin{proof}
As in the proof of Theorem \ref{thm 1}, applying the chain rule  $d|u^k|^2=2u^kdu^k$ for each $k$, we have
$$
|u^k(t)|^2=|u^k_0|^2+\int^t_0 2u^k(a^{ij}_{kr}u^r_{x^ix^j}+f^k)\,ds
$$
where the summations on $i,j,r$ are understood. Denote $c=\theta-d$.
For each  $k$, we have
\begin{eqnarray}
0&\leq& \int_{\bR^d_+}|u^k(T,x)|^2(x^1)^c
dx\nonumber\\
&=&\int_{\bR^d_+}|u^k(0,x)|^2(x^1)^cdx\nonumber\\
&& + 2\int^T_0\int_{\bR^d_+} a^{ij}_{kr}u^ku^r_{x^ix^j}(x^1)^c
dxds+2\int^T_0\int_{\bR^d_+} (M^{-1}u^k)(Mf^k)(x^1)^c
dxds .\label{2009.06.02 04:27 PM}
\end{eqnarray}

Note that, by integration by parts, the second term in
(\ref{2009.06.02 04:27 PM}) is
\begin{eqnarray}
\int^T_0\int_{\bR^d_+}\left[-2a^{ij}_{kr}u^k_{x^i}u^r_{x^j}-2
c(a^{1j}_{kr}u^r_{x^j})(M^{-1}u^k)\right](x^1)^cdxds \label{2009.09.02.6:18PM}
\end{eqnarray}
$$
\leq \int^T_0\int_{\bR^d_+}-2a^{ij}_{kr}u^k_{x^i}u^r_{x^j}(x^1)^c\,dxds+|c|\left(\kappa\|u_x\|^2_{\bL_{2,\theta}(T)}+
K^2\kappa^{-1}\|M^{-1}u\|^2_{\bL_{2,\theta}(T)}\right),
$$
for each $\kappa>0$, because for any vectors $v,w\in \bR^n$ and $\kappa>0$,
$$
|<A^{1j}v,w>|\leq |A^{1j}v||w|\leq K^j|v||w|\leq
\frac12(\kappa|v|^2+\kappa^{-1}(K^j)^2|w|^2).
$$

By summing up the terms in (\ref{2009.06.02 04:27 PM}) over $k$ and
rearranging the terms, we get
\begin{eqnarray}
&&2\int^T_0\int_{\bR^d_+}u^*_{x^i}A^{ij}u_{x^j}\;(x^1)^c
dxds\nonumber\\
&\leq&|c|\left(\kappa\|u_x\|^2_{\bL_{2,\theta}(T)}+
K^2\kappa^{-1}\|M^{-1}u\|^2_{\bL_{2,\theta}(T)}\right)
+\varepsilon \|M^{-1}u\|^2_{\bL_{2,\theta}(T)}\\
&+&c(\varepsilon)\|Mf\|^2_{\bL_{2,\theta}(T)}
+\|u(0)\|^2_{U^1_{2,\theta}},\label{eqn
2}
\end{eqnarray}
where
$\kappa,\varepsilon>0$ will be decided below. Assumption
\ref{main assumptions}(i),  inequality (\ref{eqn 2}) and the inequality
\begin{equation}
                    \label{eqn 3}
\|M^{-1}u\|^2_{L_{2,\theta}}\leq
\frac{4}{(d+1-\theta)^2}\|u_x\|^2_{L_{2,\theta}}
\end{equation}
(see Corollary 6.2 in \cite{kr99}) lead us to
 \begin{eqnarray}
 &&2\delta\|u_x\|^2_{\bL_{2,\theta}(T)}-
 |c|\left(\kappa
 +\frac{4K^2}{\kappa(d+1-\theta)^2}\right)\|u_x\|^2_{\bL_{2,\theta}(T)}\nonumber\\
&\le& c\varepsilon\|u_x\|^2_{\bL_{2,\theta}(T)}
+c(\varepsilon)\|Mf\|^2_{\bL_{2,\theta}(T)}+\|u(0)\|^2_{U^1_{2,\theta}}.\nonumber
\end{eqnarray}

Now it is enough to  take $\kappa=2K/(d+1-\theta)$ and observe that
(\ref{theta 1}) is equivalent to the condition
$$
2\delta-|c|\left(\kappa
 +\frac{4K}{\kappa(d+1-\theta)^2}\right)=2\delta-\frac{4|c|K}{d+1-\theta}>0.
 $$
 Choosing a small $\varepsilon=\varepsilon(d,\delta,\theta,K,L)$, the lemma is proved.
\end{proof}

\begin{lemma}
                     \label{a priori 2}
Let $a^{ij}_{kr}=a^{ij}_{kr}(t)$. Suppose either
\begin{equation}
             \label{con 1}
\theta\in (d-1,d], \quad
2\delta(d+1-\theta)^2-2(d+1-\theta)(d-\theta)\beta-4(d-\theta)(d+1-\theta)K^1>0
\end{equation}
or
\begin{equation}
                      \label{con 3}
\theta\in (d-1,d], \quad (\delta-\bar{\alpha})-\frac{(d-\theta)}{(d+1-\theta)}(2\delta-\beta-2\alpha)>0;
\end{equation}
or
\begin{equation}
                   \label{con 2}
\theta\in [d,d+1), \quad 8(d+1-\theta)\delta^2-(\theta-d)\beta^2>0.
\end{equation}
Let $u\in M\bH^{1}_{2,\theta}(T)$ be a solution of $($\ref{eqn
system}$)$ so that $u\in C([0,T],C^2_0((1/N,N)\times
\{x':|x'|<N\}))$ for some $N>0$. Then the assertion of Lemma \ref{a
priori 1} holds.
\end{lemma}
\begin{proof}

1. Denote $S^{1j}=(s^{1j}_{kr})=\frac12(A^{1j}+(A^{1j})^*)$ as the
symmetric part of $A^{1j}$. Then $A^{1j}=S^{1j}+\frac12 H^{1j}$, and
for any $\xi\in \bR^{d_1}$ we notice that
$\xi^*A^{1j}\xi=\xi^*S^{1j}\xi$. Let $c:=\theta-d$. Note that, by
integration by parts,
$$
\int_{\bR^d_+}u^*S^{11}u_{x^1}(x^1)^{c-1}dx=
-\frac{c-1}{2}\int_{\bR^d_+}u^* S^{11}u
(x^1)^{c-2}dx=-\frac{c-1}{2}\int_{\bR^d_+}u^* A^{11}u (x^1)^{c-2}dx
$$
and hence
\begin{eqnarray*}
 -2c\int_{\bR^d_+}u^*A^{11}u_{x^1}(x^1)^{c-1}dx&=&-2c\int_{\bR^d_+}u^*S^{11}u_{x^1}(x^1)^{c-1}dx
-c\int_{\bR^d_+}u^*H^{11}u_{x^1}(x^1)^{c-1}dx\\
&=&c(c-1)\int_{\bR^d_+}u^* A^{11}u (x^1)^{c-2}dx
-c\int_{\bR^d_+}u^*H^{11}u_{x^1}(x^1)^{c-1}dx.
 \end{eqnarray*}
Moreover, another usage of integration by parts gives us
\begin{eqnarray}
\int_{\bR^d_+}u^*S^{1j}u_{x^j}(x^1)^{c-1}dx=
-\int_{\bR^d_+}u_{x^j}^* S^{1j}u (x^1)^{c-1}dx=-\int_{\bR^d_+}
u^*(S^{1j})^*u_{x^j} (x^1)^{c-1}dx\nonumber
\end{eqnarray}
for $j\ne 1$, meaning that
$\int_{\bR^d_+}u^*S^{1j}u_{x^j}(x^1)^{c-1}dx=0$ and
\begin{eqnarray}
 -2c\int_{\bR^d_+}u^*A^{1j}u_{x^j}(x^1)^{c-1}dx=-c\int_{\bR^d_+}u^*H^{1j}u_{x^j}(x^1)^{c-1}dx.\nonumber
\end{eqnarray}
We gather the above terms to get
\begin{eqnarray}
-2c\int_{\bR^d_+}(a^{1j}_{kr}u^r_{x^j})u^k(x^1)^{c-1}dx
=c(c-1)\int_{\bR^d_+}u^* A^{11}u (x^1)^{c-2}dx
-c\int_{\bR^d_+}u^*H^{1j}u_{x^j}(x^1)^{c-1}dx,\nonumber
\end{eqnarray}
where the summation on $j$ includes $j=1$.

Now, as in the proof of Lemma \ref{a priori 1}, we have
\begin{eqnarray}
&&2\delta\|u_x\|^2_{\bL_{2,\theta}(T)}\nonumber\\
&\le&2\int^T_0\int_{\bR^d_+}u^*_{x^i}A^{ij}u_{x^j}\;(x^1)^c
dxds\nonumber\\
&\leq&\int_{\bR^d_+}|u^k(0,x)|^2x^cdx\nonumber\\
&+&c(c-1)\int^T_0\int_{\bR^d_+}a^{11}_{kr}(M^{-1}u^k)(M^{-1}u^r)(x^1)^{c}dxds
-c\int^T_0\int_{\bR^d_+}(h^{1j}_{kr}u^r_{x^j})(M^{-1}u^k) (x^1)^{c} dxds\nonumber\\
&+&2\int^T_0\int_{\bR^d_+} (M^{-1}u^k)(Mf^k)(x^1)^c
dxds.
\label{2009.06.03 06:21 PM}
\end{eqnarray}
Note that  the first and  last terms  in
the right hand side of (\ref{2009.06.03 06:21 PM}) are bounded by
$$
\varepsilon \|M^{-1}u\|^2_{\bL_{2,\theta}(T)} +c(\varepsilon)\|Mf\|^2_{\bL_{2,\theta}(T)}
+\|u(0)\|^2_{U^1_{2,\theta}}.
$$

2. If $c(c-1)\geq 0$, hence $\theta\in (d-1,d]$, then
\begin{eqnarray*}
&&c(c-1)\int^T_0\int_{\bR^d_+}a^{11}_{kr}(M^{-1}u^k)(M^{-1}u^r)(x^1)^{c}dxds\\
&\leq& c(c-1)K^1\|M^{-1}u\|^2_{\bL_{2,\theta}(T)} \leq
\frac{4}{(d+1-\theta)^2}c(c-1)K^1\|u_x\|^2_{\bL_{2,\theta}(T)}.
\end{eqnarray*}
Also,
\begin{eqnarray*}
\left|-c\int^T_0\int_{\bR^d_+}(h^{1j}_{kr}u^r_{x^j})(M^{-1}u^k)
(x^1)^{c} dxds\right|&\leq&
\frac{1}{2}|c|\left(\kappa\|u_x\|^2_{\bL_{2,\theta}(T)}
+\kappa^{-1}\beta^2\|M^{-1}u\|^2_{\bL_{2,\theta}(T)}\right)\\
&\leq& \frac{1}{2}|c|\left(\kappa
 +\frac{4\beta^2}{\kappa(d+1-\theta)^2}\right)\|u_x\|^2_{\bL_{2,\theta}(T)}
\end{eqnarray*}
for any $\kappa>0$. To minimize this we take  $\kappa=2\beta/(d+1-\theta)$, then
\begin{equation}
                      \label{eqn 1.28.1}
\left|-c\int^T_0\int_{\bR^d_+}(h^{1j}_{kr}u^r_{x^j})(M^{-1}u^k)
(x^1)^{c} dxds\right| \leq \frac{2\beta(d-\theta)}{(d+1-\theta)}\|u_x\|^2_{\bL_{2,\theta}(T)}.
\end{equation}

Thus we deduce
$$
 \left(2\delta-\frac{2\beta(d-\theta)}{(d+1-\theta)}-
 \frac{4}{(d+1-\theta)^2}c(c-1)K^1\right)\|u_x\|^2_{\bL_{2,\theta}(T)}
\leq c\varepsilon \|u_x\|^2_{\bL_{2,\theta}(T)}
+c(\varepsilon)\|Mf\|^2_{\bL_{2,\theta}(T)}+\|u(0)\|^2_{U^1_{2,\theta}}.
$$
This and  (\ref{eqn 3}) yield   a priori
 (\ref{eqn main}), since (\ref{con 1}) is equivalent to
 $$
 2\delta-\frac{2\beta(d-\theta)}{(d+1-\theta)}-
 \frac{4}{(d+1-\theta)^2}c(c-1)K^1 >0.
 $$

3. Again assume $c(c-1)\geq 0$. By (\ref{2009.06.03 06:21 PM}) and (\ref{eqn 1.28.1}),
\begin{eqnarray*}
&&2\int^T_0\int_{\bR^d_+}u^*_{x^i}\left(S^{ij}_d+S^{ij}_0\right)u_{x^j}\;(x^1)^c
dxds\\
&\leq&\int_{\bR^d_+}|u^k(0,x)|^2x^cdx\\
&+&c(c-1)\int^T_0\int_{\bR^d_+}\left(s^{11}_{d,kr}+s^{11}_{s,kr}\right)(M^{-1}u^k)(M^{-1}u^r)(x^1)^{c}dxds\\
&+&\frac{2\beta(d-\theta)}{(d+1-\theta)}\|u_x\|^2_{\bL_{2,\theta}(T)}
+\varepsilon\|M^{-1}u\|^2_{\bL_{2,\theta}(T)}+c\|Mf\|^2_{\bL_{2,\theta}(T)}.
\end{eqnarray*}

By Corollary 6.2 of \cite{kr99}, for each $t$,
$$
c(c-1)
\int s^{11}_{d,kr}(M^{-1}u^k)(M^{-1}u^r)(x^1)^{c}\,dx\leq \frac{4(d-\theta)}{(d+1-\theta)}\int_{\bR^d_+}u^*_{x^i} S^{ij}_{d}u_{x^j}\;(x^1)^c\,dx.
$$
By assumptions,
$$
2\int_{\bR^d_+} u^*_{x^i} S^{ij}_{o}u_{x^j}\;(x^1)^c\,dx\leq 2\bar{\alpha}\int_{\bR^d_+}|u_x|^2\;(x^1)^c\,dx,
$$
\begin{eqnarray*}
c(c-1)\int_{\bR^d_+} s^{11}_{0,kr}|M^{-1}u^k||M^{-1}u^r|(x^1)^c\,dx &\leq& \alpha c(c-1)\int_{\bR^d_+}|M^{-1}u|^2(x^1)^c\,dx \\
&\leq &\frac{4\alpha(d-\theta)}{(d+1-\theta)}\int_{\bR^d_+}|u_x|^2\;(x^1)^c\,dx.
\end{eqnarray*}
It follows
$$
\left[(\delta-\bar{\alpha})-\frac{(d-\theta)}{(d+1-\theta)}(2\delta-\beta-2\alpha)\right]
\|u_x\|^2_{\bL_{2,\theta}(T)}\leq \varepsilon \|u_x\|^2_{\bL_{2,\theta}(T)}+c\|Mf\|^2_{\bL_{2,\theta}(T)}
+\|u_0\|^2_{U^1_{2,\theta}}.
$$
This, (\ref{con 3}) and (\ref{eqn 3}) lead to the a priori estimate.

4. If $c(c-1)\leq 0$, hence $\theta\in [d,d+1)$, then
$$
c(c-1)\int^T_0\int_{\bR^d_+}a^{11}_{kr}(M^{-1}u^k)(M^{-1}u^r)(x^1)^{c}dxds\le
\delta c(c-1)\|M^{-1}u\|^2_{\bL_{2,\theta}(T)};
$$
for this we consider a $d_1\times d$ matrix consisting of $M^{-1}u$
as the first column and zeros for the rest, and apply the assumption
\ref{assumption 1}. Next, as before, we have
\begin{eqnarray*}
\left|-c\int^T_0\int_{\bR^d_+}(h^{1j}_{kr}u^r_{x^j})(M^{-1}u^k)
(x^1)^{c} dxds\right|&\leq&
\frac{1}{2}c\left(\kappa\|u_x\|^2_{\bL_{2,\theta}(T)}
+\kappa^{-1}\beta^2\|M^{-1}u\|^2_{\bL_{2,\theta}(T)}\right)
\end{eqnarray*}
and hence from (\ref{2009.06.03 06:21 PM}) it follows
\begin{eqnarray}
 &&2\delta\|u_x\|^2_{\bL_{2,\theta}(T)}-
 \frac12 c\;\left(\kappa\|u_x\|^2_{\bL_{2,\theta}(T)}
+\kappa^{-1}\beta^2\|M^{-1}u\|^2_{\bL_{2,\theta}(T)}\right)-\delta c(c-1)\|M^{-1}u\|^2_{\bL_{2,\theta}(T)}\nonumber\\
&\le&{\varepsilon}\|u_x\|^2_{\bL_{2,\theta}(T)}+
+c(\varepsilon) \|Mf\|^2_{\bL_{2,\theta}(T)}+\|u(0)\|^2_{U^1_{2,\theta}}.\label{2009.06.03
07:56 PM}
\end{eqnarray}
As we take
$$
\kappa=\frac{\beta^2}{2\delta(1-c)},$$ the terms with
$\|M^{-1}u\|^2_{\bL_{2,\theta}(T)}$ in the left hand side of
(\ref{2009.06.03 07:56 PM}) are canceled. Now, (\ref{con 2}) which
is equivalent to $2\delta-\frac{c\beta^2 }{4\delta(1-c)}>0$ gives us
a priori estimate (\ref{eqn main}). The lemma is proved.
\end{proof}

\begin{theorem}
                          \label{theorem half-constant}
Let $\gamma \geq 0$ and $a^{ij}_{kr}=a^{ij}_{kr}(t)$. Assume that one of (\ref{theta 1}), (\ref{con 1}), (\ref{con 3}) and (\ref{con 2}) holds.
Then for any $f\in M^{-1}\bH^{\gamma}_{2,\theta}(T)$ and $ u_0\in
U^{\gamma+2}_{2,\theta}$, system (\ref{eqn system}) admits a unique
solution $u\in \frH^{\gamma+2}_{2,\theta}(T)$, and for this solution
\begin{equation}
                        \label{a priori}
\|M^{-1}u\|_{\bH^{\gamma+2}_{2,\theta}(T)}\leq
c\|Mf\|_{\bH^{\gamma}_{2,\theta}(T)}+c\|u_0\|_{U^{\gamma+2}_{2,\theta}},
\end{equation}
where $c=c(d,\delta,\theta,K,L)$.
\end{theorem}
\begin{proof}
1. By Theorem 3.3 in \cite{KL2}, for each $k$, the equation
$$
u^k_t=\delta \Delta u^k+f^k, \quad u^k(0)=u^k_0
$$
has a solution $u^k\in \frH^{\gamma+2}_{2,\theta}(T)$. As in the
proof of Theorem \ref{thm 1} we only need to show that estimate
(\ref{a priori}) holds given that a solution already exists.

2. By Theorem 2.9 in \cite{KL2}, for any nonnegative integer $n\geq
\gamma+2$, the set
$$
\frH^{n}_{2,\theta}(T) \cap
\bigcup_{N=1}^{\infty}C([0,T],C^n_0((1/N,N)\times
\{x':|x'|<N\}))
$$
is everywhere dense in $\frH^{\gamma+2}_{p,\theta}(T)$ and we may
assume that $u$ is sufficiently smooth in $x$ and vanishes near the
boundary. Thus  a priori estimate (\ref{a
priori}) follows from  Lemma \ref{lemma 2.05.10}, Lemma \ref{a priori 1} and Lemma \ref{a priori
2}. The theorem is proved.

\end{proof}

Here is the main result of this section.

\begin{theorem}
                 \label{theorem half}
Let $\gamma \geq 0$ and Assumption \ref{assumption theta} hold. Assume that  for each $t$
$$
|a^{ij}_{kr}(t,\cdot)|^{(0)*}_{\gamma_+}
+|b^{i}_{kr}(t,\cdot)|^{(1)*}_{\gamma_+}
+|c_{kr}(t,\cdot)|^{(2)*}_{\gamma_+}  \leq L
$$
and
$$
|a^{ij}_{kr}(t,x)-a^{ij}_{kr}(t,y)| +
|Mb^i_{kr}(t,x)|+|M^2c_{kr}(t,x)|<\kappa
$$
for all $x,y\in \bR^d_+$ with $|x-y|\leq x^1\wedge y^1$. Then there
exists $\kappa_0=\kappa_0(d,\theta,\delta,K,L)$ so that if $\kappa\leq
\kappa_0$, then for any $f\in M^{-1}\bH^{\gamma}_{2,\theta}(T)$,
 and $u_0\in
U^{\gamma+2}_{2,\theta}$, system (\ref{eqn system2}) admits a unique solution $u\in
\frH^{\gamma+2}_{2,\theta}(T)$, and furthermore
\begin{equation}
                                  \label{eqn 9.9.2}
 \|u\|_{\frH^{\gamma+2}_{2,\theta}(T)}\leq
 c\|Mf\|_{\bH^{\gamma}_{2,\theta}(T)}
 +c\|u_0\|_{U^{\gamma+2}_{2,\theta}}
 \end{equation}
 where $c=c(d,d_1,\delta,\theta,K,L)$.
 \end{theorem}

 To prove Theorem \ref{theorem half} we use
the following lemmas  taken from \cite{KK2}.

\begin{lemma}
                                         \label{lemma 8.26.10}
Let constants $C,\delta\in(0,\infty)$,   a function
 $u\in H^{\gamma}_{p,\theta} $, and $q$ be the smallest integer such
that $|\gamma|+2\leq q$.

(i) Let $\eta_{n}\in C^{\infty}(\bR^{d}_{+})$, $n=1,2,...$, satisfy
\begin{equation}
                                                   \label{8.26.11}
\sum_{n}M^{|\alpha|} |D^{\alpha} \eta_n |\leq C
\quad\text{in}\quad\bR^{d}_{+}
\end{equation}
for    any multi-index $\alpha$ such that $ 0\leq |\alpha| \leq q$.
Then
$$
\sum_{n} \|\eta_{n}u\|^{p}_{H^{\gamma}_{p,\theta} } \leq
NC^{p}\|u\|^{p}_{H^{\gamma}_{p,\theta} },
$$
where the constant $N$ is independent of $u$, $\theta$, and $C$.

(ii) If in addition to the condition in (i)
\begin{equation}
                                                   \label{1.5.2}
\sum_{n} \eta_{n} ^{2}\geq\delta\quad\text{on}\quad\bR^{d}_{+},
\end{equation}
then
\begin{equation}
                                                   \label{11.25.1}
\|u\|^{p}_{H^{\gamma}_{p,\theta} }\leq N\sum_{n}
\|\eta_{n}u\|^{p}_{H^{\gamma}_{p,\theta} },
\end{equation}
where the constant $N$ is independent of $u$ and $\theta$.
\end{lemma}

The reason the first inequality in (\ref{11.14.1}) below is written
for $\eta_n^4$ (not for $\eta_n^2$) as in the above lemma is to have
the possibility to apply Lemma \ref{lemma 8.26.10} to $\eta_n^2$.
Also observe that obviously $\sum a^2 \leq (\sum |a|)^2$.

\begin{lemma}
                                            \label{lemma 11.14.1}
For each $\varepsilon>0$ and $q=1,2,...$ there exist non-negative
functions $\eta_{n}\in C^{\infty}_{0}(\bR^{d}_{+})$, $n=1,2,...$
such that (i) on $\bR^{d}_{+}$ for each multi-index $\alpha$ with
$1\leq|\alpha|\leq q$ we have
\begin{equation}
                                             \label{11.14.1}
\sum_{n}\eta^{4}_{n}\geq1,\quad \sum_{n} \eta _{n} \leq N(d),
\quad\sum_{n}M^{|\alpha|}|D^{\alpha}\eta_{n}|\leq\varepsilon;
\end{equation}

(ii) for any $n$ and $x,y\in\text{\rm supp}\,\eta_{n}$ we have $
|x-y|\leq N ( x^{1}\wedge y^{1})$, where
$N=N(d,q,\varepsilon)\in[1,\infty)$.

\end{lemma}

\begin{lemma}
                                           \label{lemma 1.2.2}
Let $p\in(1,\infty)$, $\gamma,\theta\in\bR$. Then there exists a
constant $N=N(\gamma,|\gamma|_+,p,d)$   such that  if $f\in
H^{\gamma}_{p,\theta}$ and $a$ is a function with finite norm
$|a|^{(0)*}_{|\gamma|_+,\bR^{d}_{+}}$, then
 \begin{equation}
                                                 \label{8.19.05}
\|af\|_{H^{\gamma}_{p,\theta}} \leq N
|a|^{(0)*}_{|\gamma|_+,\bR^{d}_{+}} \|f\|_{H^{\gamma }_{p,\theta}}.
\end{equation}
In addition,

(i) if $\gamma=0,1,2,...$, then
 \begin{equation}
                                            \label{1.24.06}
\|af\|_{H^{\gamma}_{p,\theta}} \leq N \sup_{\bR^{d}_{+}}|a|\,
\|f\|_{H^{\gamma }_{p,\theta}}+ N_0\|f\|_{H^{\gamma-1}_{p,\theta}}
\sup_{\bR^{d}_{+}}\sup_{1\leq|\alpha|\leq\gamma}
|M^{|\alpha|}D^{\alpha}a |,
\end{equation}
where $N_0=0$ if $\gamma=0$, and $N_0=N_0(\gamma,d)>0$ otherwise.

(ii) if $\gamma$ is not integer, then
\begin{equation}
                                            \label{1.24.07}
\|af\|_{H^{\gamma}_{p,\theta}} \leq N  (\sup_{\bR^{d}_{+}}|a|)^{s}
(|a|^{(0)*}_{|\gamma|_+})^{1-s} \|f\|_{H^{\gamma }_{p,\theta}},
\end{equation}
where $s:=1-\frac{|\gamma|}{{|\gamma|_+}} > 0$.

\end{lemma}

{\bf{Proof of Theorem \ref{theorem half}}}

 We closely follow the proof of Theorem 2.16 of \cite{KK}.
 As usual, for simplicity, we assume $u_0 =0$.
Also  having  the method of continuity in mind, we convince
ourselves  that
 to prove the theorem  it  suffices to
show that there exist  $\kappa_{0}$  such that
 the a priori estimate (\ref{eqn 9.9.2}) holds
 given that the solution already exists and $\kappa
\leq\kappa_{0}$.
  We divide the
proof into two cases. This is because    if $\gamma$ is an
 integer we use (\ref{1.24.06}), and otherwise  we  use (\ref{1.24.07}).

 {\bf Case 1}: $\gamma=0$  or $\gamma$ is not integer.
 Take the least integer $q\geq|\gamma|+4$. Also take an
$\varepsilon\in(0,1)$ to be specified later and take a sequence of
functions $\eta_{n}$, $n=1,2,...$, from Lemma \ref{lemma 11.14.1}
corresponding to $\varepsilon,q$.
 Then by Lemma \ref{lemma 8.26.10}, we have
\begin{equation}
                                                \label{8.28.15}
\|M^{-1}u\|_{\bH^{\gamma+2}_{2,\theta}(T)}^{2} \leq
N\sum_{n=1}^{\infty}
\|M^{-1}u\eta^{2}_{n}\|_{\bH^{\gamma+2}_{2,\theta}(T)}^{2}.
\end{equation}
 For any $n$ let  $x_{n}$ be a point in $\text{supp}\,\eta_{n}$
and $a^{ij}_{n,kr}(t)=a^{ij}_{kr}(t,x_{n})$.
 From  (\ref{eqn system2}), we easily have
$$
 (u^k\eta^{2}_{n})_t=
a^{ij}_{n,kr}(u^r\eta^{2}_{n})_{ x^ix^j}+M^{-1}f^k_{n},
$$
where
$$
f^k_{n}=(a^{ij}_{kr}-a^{ij}_{n,kr}) \eta^{2}_{n} Mu^r_{x^ix^j}
 -2a^{ij}_{n,kr}M(\eta^{2}_{n})_{ x^i}u^r_{x^j}
-a^{ij}_{n,kr}M^{-1}u^rM^{2}(\eta^{2}_{n})_{ x^ix^j}
$$
$$
+\eta_{n}^{2}Mb^{i}_{kr}u^r_{x^{i}}
+\eta_{n}^{2}M^{2}c_{kr}M^{-1}u^r +Mf^k\eta^{2}_{n}.
$$

By Theorem \ref{theorem half-constant}, for each n,
\begin{equation}
                                                 \label{8.28.20}
\|M^{-1}u\eta^{2}_{n}\|_{\bH^{\gamma+2}_{2,\theta}(T)}^{2} \leq N
\|f_{n}\|^{2}_{\bH^{\gamma}_{2,\theta}(T)}
\end{equation}
and  by (\ref{1.24.07}),
\begin{equation}
                                                             \label{1.24.01}
\|(a^{ij}_{kr}-a^{ij}_{n,kr}) \eta^{2}_{n} Mu_{x^ix^j}\|
_{\bH^{\gamma}_{p,\theta}} \leq N\| \eta _{n} Mu_{xx}\|
_{\bH^{\gamma}_{p,\theta}} \sup_{t,x}|(a^{ij}_{kr}-
a^{ij}_{n,kr})\eta_{n}|^{s},
\end{equation}
 where $s=1$ if $\gamma=0$, and $s=1-\frac{\gamma}{\gamma_+} >0$ otherwise.

By Lemma \ref{lemma 11.14.1}(ii), for each   $n$ and
$x,y\in\text{supp}\,\eta_{n}$
  we have
$|x-y|\leq N(\varepsilon)(x^{1}\wedge y^{1})$, where
$N(\varepsilon)=N(d,q,\varepsilon)$, and we can easily find not more
than $N(\varepsilon)+2 \leq 3N(\varepsilon)$ points $x_i$ lying on
the straight segment
 connecting $x$ and $y$ and including $x$ and $y$, such that $|x_i-x_{i+1}|\leq
x^1_{i} \wedge x^1_{i+1}$. It follows  from our assumptions
$$
\sup_{\omega,t,x}|(a^{ij}_{kr}-a^{ij}_{n,kr})\eta_{n}|
 \leq 3N(\varepsilon)\kappa.
$$
We substitute this  to (\ref{1.24.01})  and get
$$
\|(a^{ij}_{kr}-a^{ij}_{k,rn}) \eta^{2}_{n}
Mu^r_{x^ix^j}\|_{\bH^{\gamma}_{2,\theta}(T)} \leq N
N(\varepsilon)\kappa^{s}\| \eta _{n}
Mu_{xx}\|_{\bH^{\gamma}_{2,\theta}(T)}.
$$
Similarly,
$$
\| \eta^{2}_{n} Mb^{i}_{kr}u^r_{x^i}\|
_{\bH^{\gamma}_{2,\theta}(T)}+ \| \eta^{2}_{n}
M^{2}c_{kr}M^{-1}u^r\| _{\bH^{\gamma}_{2,\theta}(T)} \leq N
N(\varepsilon)
\kappa^{s}(\|\eta_{n}u_{x}\|_{\bH^{\gamma}_{2,\theta}(T)}
+\|\eta_{n}M^{-1}u\| _{\bH^{\gamma}_{2,\theta}(T)}).
$$

Coming back to (\ref{8.28.20}) and (\ref{8.28.15}) and using Lemma
\ref{lemma 8.26.10}, we conclude
$$
\|M^{-1}u\|_{\bH^{\gamma+2}_{2,\theta}(T)}^{2} \leq NN(\varepsilon)
\kappa^{2s}(\|M u_{xx}\|_{\bH^{\gamma}_{2,\theta}(T)}^{2} +
\|u_{x}\|_{\bH^{\gamma}_{2,\theta}(T)}^{2} +
\|M^{-1}u\|_{\bH^{\gamma}_{2,\theta}(T)}^{2})
$$
\begin{equation}
                                               \label{11.19.2}
+NC^{2}\left(\|u_{x}\|_{\bH^{\gamma }_{2,\theta}(T)}^{2}
+\|M^{-1}u\|_{\bH^{\gamma+1}_{2,\theta}}^{2}\right)+N
\|Mf\|_{\bH^{\gamma }_{2,\theta}}^{2},
\end{equation}
where
$$
C=\sup_{\bR^{d}_{+}}\sup_{|\alpha|\leq q-2}
\sum_{n=1}^{\infty}M^{|\alpha|}(|D^{\alpha}(M(\eta_{n}^{2})_{x})|
+|D^{\alpha}(M^{2}(\eta_{n}^{2})_{xx})|).
$$
By construction, we have $C\leq N\varepsilon$. Furthermore (see, Lemma \ref{collection half})
\begin{equation}
                                               \label{11.19.1}
\|u_{x}\|_{H^{\gamma+1}_{2,\theta}} \leq
N\|M^{-1}u\|_{H^{\gamma+2}_{2,\theta}},\quad
\|Mu_{xx}\|_{H^{\gamma}_{2,\theta}} \leq
N\|M^{-1}u\|_{H^{\gamma+2}_{2,\theta}}.
\end{equation}

Hence (\ref{11.19.2}) yields
$$
\|M^{-1}u\|_{\bH^{\gamma+2}_{2,\theta}(T)}^{2} \leq
N_{1}(N(\varepsilon)\kappa^{2s}+\varepsilon^{2})
\|M^{-1}u\|_{\bH^{\gamma+2}_{2,\theta}(T)}^{2} +N ( \|Mf\|_{\bH^{\gamma
}_{2,\theta}(T)}^{2}).
$$
Finally  to get the a priori estimate,  it's enough to
choose first
 $\varepsilon$ and then $\kappa_{0}$,
so that $N_{1}(N(\varepsilon)\kappa^{2s} +\varepsilon^{2})\leq1/2$
for $\kappa\leq\kappa_{0}$.

 {\bf Case 2}: $\gamma \in\{1,2,...\}$.
Proceed as in Case 1 with $\varepsilon=1$ and
arrive at (\ref{8.28.20}) which is
$$
\|M^{-1}u\eta^{2}_{n}\|_{\bH^{\gamma+2}_{2,\theta}(T)}^{2} \leq
N\|f_{n}\|^{2}_{\bL_{2,\theta}(T)}.
$$
Now we use (\ref{1.24.06}) to get
$$
\|(a^{ij}_{kr}-a^{ij}_{k,rn}) \eta^{2}_{n}
Mu^r_{x^ix^j}\|_{\bH^{\gamma}_{2,\theta}(T)} \leq N
\kappa\| \eta _{n}
Mu_{xx}\|_{\bH^{\gamma}_{2,\theta}(T)}+N\| \eta _{n}
Mu_{xx}\|_{\bH^{\gamma-1}_{2,\theta}(T)}.
$$

From this point by following the arguments in case 1, one easily gets
 \begin{equation}
                      \label{eqn 2.06.2}
\|M^{-1}u\|_{\bH^{\gamma+2}_{2,\theta}(T)} \leq N_{1} \kappa
\|M^{-1}u\|_{\bH^{\gamma+2}_{2,\theta}(T)}
 + N_2 \|M^{-1}u\|_{\bH^{\gamma+1}_{2,\theta}(T)} +
N\|Mf\|_{\bH^{\gamma}_{2,\theta}(T)}.
\end{equation}
This and the embedding inequality
$$
\|M^{-1}u\|_{H^{\gamma+1}_{2,\theta}}\leq \frac{1}{2N_2} \|M^{-1}u\|_{H^{\gamma+2}_{2,\theta}}+
N(N_2,\gamma)\|M^{-1}u\|_{H^2_{2,\theta}}
$$
yield
\begin{equation}
                                                            \label{1.24.03}
\|M^{-1}u\|_{\bH^{\gamma+2}_{2,\theta}(T)} \leq  2N_1 \kappa \|M^{-1}u\|_{\bH^{\gamma+2}_{2,\theta}(T)}+N
\|M^{-1}u\|_{\bH^2_{2,\theta}(T)} + N\|Mf\|_{\bH^{\gamma}_{2,\theta}(T)}.
\end{equation}
Now take  $\kappa_0$ is from Case 1 for $\gamma=0$, then it is enough to  assume
$\kappa \leq \kappa_0\wedge 1/{(4N_1)}$, because
by the result of Case 1,
$$
\|M^{-1}u\|_{\bH^2_{2,\theta}(T)}\leq N\|Mf\|_{\bL_{2,\theta}(T)}.
$$
 The theorem is proved.

\mysection{Proof of Theorem \ref{main theorem on domain}}
                              \label{section 5}
By Theorem 2.10 in \cite{KK2}, for each $k$, $f^k \in \psi^{-1}\bH^{\gamma}_{2,\theta}(\cO,T)$ and
$u^k_0\in U^{\gamma+2}_{2,\theta}(\cO)$, the equation
$$
u^k_t=\Delta u^k +f^k, \quad u^k(0)=u^k_0(0)
$$
has a unique solution $u\in \frH^{\gamma+2}_{2,\theta}(T)$, and furthermore
$$
\|\psi^{-1}u^k\|_{\bH^{\gamma+2}_{2,\theta}(\cO,T)}\leq c\|\psi f^k\|_{\bH^{\gamma}_{2,\theta}(\cO,T)}
+c\|u^k_0\|_{U^{\gamma+2}_{2,\theta}(\cO)}.
$$
Thus to prove the theorem we only
need to  prove that   (\ref{a priori domain})  holds given that  a
solution $u\in \frH^{\gamma+2}_{2,\theta}(\cO,T)$ already exists. As
usual we assume  $u_0=0$.
 Let $x_0 \in \partial \cO$ and $\Psi$ be a function from
Assumption \ref{assumption domain}. In \cite{KK2} it is shown that
 $\Psi$ can be chosen in such a way that  for any non-negative integer $n$
\begin{equation}
                                                          \label{2.25.03}
|\Psi_{x}|^{(0)}_{n,B_{r_0}(x_0)\cap \cO} +
 |\Psi^{-1}_{x}|^{(0)}_{n,J_{+}} < N(n)<  \infty
\end{equation}
and
\begin{equation}
                                                     \label{2.25.02}
\rho(x)\Psi_{xx}(x) \to 0 \quad \text{as}\quad  x\in
B_{r_0}(x_0)\cap \cO,
 \text{and} \,\,\,  \rho(x) \to 0,
\end{equation}
where the constants $N(n)$ and the
 convergence in (\ref{2.25.02}) are independent of  $x_0$.

 Define $r=r_{0}/K_{0}$
and fix smooth functions $\eta \in C^{\infty}_{0}(B_r ), \varphi\in
C^{\infty}(\bR)$ such that $ 0 \leq \eta, \varphi \leq 1$, and
 $\eta=1$ in $B_{r/2} $,  $\varphi(t)=1$ for $t\leq -3$, and
$\varphi(t)=0$
 for $t\geq-1$. Observe that
$\Psi(B_{r_0}(x_0))$ contains $B_r $.
 For $n=1,2,... $, $t>0$, $x\in\bR^{d}_{+}$ introduce
$\varphi_{n}(x)=\varphi(n^{-1}\ln x^1)$,
$$
\hat{a}^{ij,n}(t,x):= \eta(x)\varphi_n(x)\left(\sum_{l,m=1}^d
  a^{lm}(t,\Psi^{-1}(x))\cdot\partial_l\Psi^{i}(\Psi^{-1}(x))\cdot\partial_m\Psi^{j}(\Psi^{-1}(x))\right)  +
\delta^{ij}(1- \eta(x)\varphi_n(x) )I,
$$
\begin{eqnarray*}
\hat{b}^{i,n}(t,x) &:=&\eta(x) \varphi_{n}(x)\Big[
\sum_{l,m}a^{lm}(t,\Psi^{-1}(x))\cdot
\partial_{lm}\Psi^{i}(\Psi^{-1}(x))+\sum_{l}b^{l}(t,\Psi^{-1}(x))\cdot\partial_l\Psi^{i}(\Psi^{-1}(x))\Big],
\end{eqnarray*}
$$
 \hat{c}^{n}(t,x) :=\eta(x) \varphi_{n}(x)c(t,\Psi^{-1}(x)).
 $$
Then  by Assumption \ref{assumption regularity}(iii) and (\ref{2.25.03}),
one can  show that there is a  constant $L'$ independent of $n$ and $x_0$ such that
$$
|\hat{a}^{ij,n}(t,\cdot)|^{(0)*}_{\gamma_+}
+|\hat{b}^{i,n}_{kr}(t,\cdot)|^{(1)*}_{\gamma_+}
+|\hat{c,n}(t,\cdot)|^{(2)*}_{\gamma_+}  \leq L'.
$$

Take  $\kappa_0$ from Theorem \ref{theorem
half} corresponding to $d,d_1,\theta,\delta, K$ and $L'$.
 Observe that $\varphi_n(x)=0$ for $x^1 \geq e^{-n}$. Also
 (\ref{2.25.02}) implies $x^{1}\Psi_{xx}(\Psi^{-1}(x)) \to 0$ as $x^1\to 0$.
 Using these facts and
Assumption \ref{assumption regularity}(ii), one can find
 $n>0$ independent of $x_0$ such that
$$
|\hat{a}^{ij,n}_{kr}(t,x)-\hat{a}^{ij,n}_{kr}(t,y)|
+ (x^1)|\hat{b}^{i,}_{kr}(t,x)|+ x^2|\hat{c}^n_{kr}(t,x)| \leq \kappa_0,
$$
whenever $t>0, x,y\in \bR^d_+$ and $|x-y|\leq x^1 \wedge y^1$.
  Now we fix  a $\rho_0  <r_0  $ such that
$$
\Psi(B_{\rho_0}(x_0)) \subset B_{r/2} \cap \{x:x^1 \leq e^{-3n }\}.
$$

 Let $\zeta$ be a smooth function
with support in $B_{\rho_0}(x_0)$ and denote
$v:=(u\zeta)(\Psi^{-1})$ and continue $v$ as zero in
$\bR^{d}_{+}\setminus\Psi(B_{\rho_0}(x_0))$. Since
$\eta\varphi_{n}=1$ on $\Psi(B_{\rho_0}(x_0))$, the function  $v$
satisfies
$$
 v^k_t = \hat{a}^{ij,n}_{kr}v^r_{x^i x^j} + \hat{b}^{i,n}_{kr}v^r_{x^i} +
\hat{c}^{n}_{kr}v^r + \hat{f}^k
$$
where
$$\hat{f}^k =\tilde{f}^k(\Psi^{-1}), \quad \tilde{f}^k=
-2a^{ij}_{kr}u^r_{x^{i}}\zeta_{x^{j}}
-a^{ij}_{kr}u^r\zeta_{x^{i}x^{j}}-b^{i}_{kr}u^r\zeta_{x^{i}} +\zeta f^k.
$$

Next we observe  that by Lemma \ref{lemma 10.3.1} and
  Theorem 3.2 in \cite{Lo2} (or see \cite{KK2})
for any $\nu,\alpha \in \bR $ and $h \in
\psi^{-\alpha}H^{\nu}_{p,\theta}(\cO)$ with support in
$B_{\rho_0}(x_0)$
\begin{equation}
                                                          \label{1.28.01}
\|\psi^{\alpha}h\|_{H^{\nu}_{p,\theta}(\cO)} \sim
\|M^{\alpha}h(\Psi^{-1})\|_{H^{\nu}_{p,\theta}}.
\end{equation}
Therefore we  conclude that
 $v\in \frH^{\gamma+2}_{2,\theta}(T)$, and  by Theorem \ref{theorem half}  we have, for any $t\leq T$,
$$
\|M^{-1}v\|_{\bH^{\gamma+2}_{2,\theta}(t)} \leq N \|M\hat{f}
\|_{\bH^{\gamma}_{2,\theta}(t)}.
$$
By using (\ref{1.28.01}) again we obtain
 \begin{eqnarray*}
\|\psi^{-1}u\zeta\|_{\bH^{\gamma+2}_{2,\theta}(\cO,t)} &\leq& N
\|a\zeta_x \psi u_x\|_{\bH^{\gamma}_{2,\theta}(\cO,t)} + N
\|a\zeta_{xx}\psi u\|_{\bH^{\gamma}_{2,\theta}(\cO,t)}\\
&+& N \|\zeta_x \psi b u\|_{\bH^{\gamma}_{2,\theta}(\cO,t)}  +  N \|\zeta \psi
f\|_{\bH^{\gamma}_{2,\theta}(\cO,t)}.
\end{eqnarray*}

Next, we easily check that
$$
|\zeta_x a(t,\cdot)|^{(0)}_{|\gamma|_+},\,\, |\zeta_{xx}\psi
a(t,\cdot)|^{(0)}_{|\gamma|_+},\,\, |\zeta_x \psi
b(t,\cdot)|^{(0)}_{|\gamma|_+}
$$
are bounded on $[0,T]$, and    conclude
$$
\|\psi^{-1}u\zeta\|_{\bH^{\gamma+2}_{2,\theta}(\cO,t)} \leq N \|\psi
u_x\|_{\bH^{\gamma}_{2,\theta}(\cO,t)} + N
\|u\|_{\bH^{\gamma}_{2,\theta}(\cO,t)}
+N \|\psi f\|_{\bH^{\gamma}_{2,\theta}(\cO,t)}.
$$

Finally, to estimate the norm
 $\|\psi^{-1} u\|_{\bH^{\gamma+2}_{2,\theta}(\cO,t)}$,
 we introduce a partition of unity $\zeta_{(i)}, i=0,1,2,...,M$ such
that $\zeta_{(0)} \in C^{\infty}_0(\cO)$ and
  $\zeta_{(i)} \in C^{\infty}_0(B_{\rho_0}(x_i))$,
$ x_i \in \partial \cO$ for $i\geq1$.
   Observe that since
$u\zeta_{(0)}$ has compact in $\cO$, we get
$$
\|\psi^{-1}u\zeta_{(0)}\|_{\bH^{\gamma+2}_{2,\theta}(\cO,t)}\sim
\|u\zeta_{(0)}\|_{\bH^{\gamma+2}_{2}(t)}.
$$
Thus we can estimate
$\|\psi^{-1} u\zeta_{(0)}\|_{\bH^{\gamma+2}_{2,\theta}(\cO,t)}$ using
Theorem \ref{thm 2} and the other norms as above.

By summing
up those estimates we get
$$
\|\psi^{-1} u\|_{\bH^{\gamma+2}_{2,\theta}(\cO,t)} \leq N
 \|\psi u_x\|_{\bH^{\gamma}_{2,\theta}(\cO,t)}+
N\|u\|_{\bH^{\gamma}_{2,\theta}(2,t)}
+ N \|\psi f\|_{\bH^{\gamma}_{2,\theta}(\cO,t)}.
$$
Furthermore, we know   that
$$
\|\psi u_x\|_{ H^{\gamma}_{2,\theta}(\cO)} \leq N \|u\|_{
H^{\gamma+1}_{2,\theta}(\cO)}.
$$
Therefore it follows
\begin{eqnarray*}
\|u\|^2_{\frH^{\gamma+2}_{2,\theta}(\cO,t)} &\leq& N
\|u\|^2_{\bH^{\gamma+1}_{2,\theta}(\cO,t)} + N \|\psi
f\|^2_{\bH^{\gamma}_{2,\theta}(\cO,t)}\\
&\leq & N \int^t_0 \|u\|^2_{\frH^{\gamma+2}_{2,\theta}(\cO,s)}\,ds+N \|\psi
f\|^2_{\bH^{\gamma}_{2,\theta}(\cO,t)}
\end{eqnarray*}
where Lemma \ref{lemma 15.05} is used for the second inequality.  Now (\ref{a priori domain})  follows from
Gronwall's inequality.  The theorem is proved.

\mysection{Proof of Theorem \ref{theorem elliptic}}
                            \label{section 6}

 Again we only show that a priori estimate (\ref{a priori elliptic}) holds
given that a solution $u\in \psi H^{\gamma+2}_{2,\theta}(\cO)$
already exists. By  (\ref{03.04.01}) it follows that  $\psi$ is a point-wise multiplier in $H^{\nu}_{p,\theta}(\cO)$ for any $\nu$ and $p$. Thus
\begin{equation}
                                \label{eqn 6.24}
\|u\|_{U^{\gamma+2}_{2,\theta}(\cO)}:=\|u\|_{H^{\gamma+1}_{2,\theta}(\cO)}\leq
c(\theta,\gamma)\|\psi^{-1}u\|_{H^{\gamma+1}_{2,\theta}(\cO)}.
\end{equation}
Note that  $v^k:=u^ke^{\lambda^k t}$ satisfies
$$
v^k_t=a^{ij}_{kr}v^r_{x^ix^j}+b^{i}_{kr}v^r_{x^i}+c_{kr}v^r+f^ke^{\lambda^k t}.
$$
By (\ref{a priori domain}) and (\ref{eqn 6.24}),
$$
g_1(T)\|\psi^{-1}u\|_{H^{\gamma+2}_{2,\theta}(\cO)}\leq
ce^{cT}\left(\|\psi^{-1}u\|_{H^{\gamma+2}_{2,\theta}(\cO)}+g_2(T)\|\psi
f\|_{H^{\gamma}_{2,\theta}(\cO)}\right),
$$
where
$$
g_1(T)=\left(\int^T_0 e^{2t\min\{\lambda^k\}}dt\right)^{1/2}, \quad
g_2(T)=\left(\int^T_0 e^{2t\max\{\lambda^k\}}dt\right)^{1/2}.
$$
If $\min\{\lambda^k\}>c$, then the ratio $ce^{cT}/{g_1(T)}$ tends to
zero as $T\to \infty$. Then after finding a $T$ such that this ratio
is less than $1/2$ one gets (\ref{a priori elliptic}). The theorem is
proved.

\end{document}